\documentclass[11pt]{article}

\usepackage[T1]{fontenc}
\usepackage[utf8]{inputenc}
\usepackage{amsmath,amssymb,amsthm,mathtools}
\usepackage{bm}
\usepackage{enumitem}
\usepackage{graphicx}
\usepackage{geometry}
\usepackage{hyperref}
\usepackage{xcolor}

\usepackage{wrapfig}
\usepackage{tikz}
\usepackage{placeins}
\usetikzlibrary{calc}

\newtheorem{theorem}{Theorem}[section]
\newtheorem{proposition}[theorem]{Proposition}
\newtheorem{lemma}[theorem]{Lemma}
\newtheorem{corollary}[theorem]{Corollary}
\theoremstyle{definition}
\newtheorem{definition}[theorem]{Definition}
\theoremstyle{remark}
\newtheorem{remark}[theorem]{Remark}

\newcommand{\R}{\mathbb{R}}

\newcommand{\Sseg}{\mathcal{S}}
\newcommand{\Hcal}{\mathcal{H}}

\newcommand{\dd}{\,\mathrm{d}}

\title{Singularities and Collision Criteria in the Gravitational Straight-Segment Problem with Arbitrary Positive Density}
\author{F. Crespo, \; E. Mart\'{\i}nez,\; J. Vidarte,\;  J.L. Zapata, \; A. Abdrabou}
\date{}

\begin{document}

\maketitle


\begin{abstract}
We study singularities and collision criteria for a massless particle around a fixed finite straight segment with an arbitrary positive continuous density. We first show that every singularity forces the trajectory to approach the material segment. Using the axial symmetry of the potential, we prove that any trajectory approaching the segment is confined to an invariant plane containing it. Moreover, we show that every singularity is either a collision singularity or an ejection singularity. Finally, in the invariant plane, we introduce an angular-cone condition for positions and momenta directed toward the segment. We prove that if the initial position-momentum pair belongs to the angular cone, then the trajectory collides with the segment in finite time, and we derive an upper bound for the collision time. We further prove that the final arc of any collision satisfies the
angular-cone condition.
\end{abstract}


\section{Introduction}
\label{sec:introduction}

In this paper, we study the motion of a massless particle attracted by a fixed finite straight segment carrying an arbitrary positive continuous density. We characterize the finite-time singularities of the associated flow, identify the geometric obstruction to collision, and derive sufficient conditions for collision in the three-dimensional problem. This analysis is motivated by both mathematical and applied considerations. From the mathematical perspective, classifying singularities determines whether a finite-time failure of continuation corresponds to a collision or to another mechanism. For an extended body, this question is more subtle than in the classical Kepler problem: the collision set is not a single point, the potential has a different local singular behavior, and angular momentum does not by itself determine the fate of the particle. From the applied perspective, sufficient collision conditions can provide analytical screening criteria for identifying potential collision courses and may inform the design of collision-avoidance maneuver algorithms for proximity operations and surface transfers near small bodies
\cite{MashikuCarpenter2016,LiuYangLi2021}.

A finite straight segment is a natural idealized model for an elongated massive body and is motivated by the study of spacecraft motion near small, elongated and irregular asteroids. In such environments, the shape and mass distribution of the body influence nearby trajectories, gravity-field modeling, and proximity operations. Straight-segment models provide a useful intermediate setting between high-fidelity representations \cite{MullerSjogren1968,WernerScheeres1996} and analytically tractable gravitational fields \cite{Bartczak2003,Bartczak2006,Riaguas1999,RiaguasElipeLara1999,ZengJiangLiBaoyin2015}. The model captures the essential geometric feature of elongation while retaining enough structure for a rigorous dynamical analysis. In particular, it has axial symmetry and therefore admits a conserved angular momentum about the segment axis, whereas its gravitational field is not central. The fixed straight-segment problem has been studied in connection with periodic orbits, stability, invariant sets, and non-integrability \cite{Duboshin1959,RiaguasElipeMoratalla2001,ArribasElipe2001,PalacianYanguasRomero2006}. More recent works have considered non-homogeneous densities, rotating segments, periodic orbits, and KAM tori \cite{Elisa2024a,Elisa2024c, MartinezVidarteRotation}. The singularity problem considered here is also connected with the analysis developed in 
\cite{MartinezThesis}, where the author focused on the case of linear density.

Singularities have also been studied for other extended mass distributions. Azevedo, Cabral, and Ontaneda investigated the fixed homogeneous circle problem and proved that all finite-time singularities are collision singularities \cite{Azevedo2003,Azevedo2007}. Alberti and Vidal studied the gravitational attraction generated by homogeneous straight segments, annular disks, and disks \cite{AlbertiVidal2010}. In the straight-segment case with constant density, their collision analysis is carried out on the equatorial plane. They also considered a related problem involving a central mass point  surrounded by a solid circular ring \cite{AlbertiVidal2019}. In their planar reductions, nonzero angular momentum produces a centrifugal barrier, and collision requires the vanishing of the corresponding planar angular momentum. This leads to a radial reduction closely related to the familiar mechanism in the Kepler
problem.

The present work addresses the segment problem without imposing a planar or equatorial restriction. We first analyze the full three-dimensional flow and only subsequently derive a planar reduction from the necessary conditions for approach to the segment. In addition, we allow an arbitrary positive continuous density rather than restricting the mass distribution to a homogeneous or prescribed profile. Thus, the invariant plane used in the later part of the paper is a consequence of the singularity analysis, rather than an assumption imposed on the initial
problem.

The paper is organized as follows. Section \ref{sec:model} introduces the gravitational potential and Hamiltonian formulation. Section \ref{sec:singularities} characterizes the potential singularity and establishes the necessary approach to the material segment at finite singular times. Section \ref{sec:CollisionConditions} derives the necessary angular-momentum condition, the invariant planar reduction, and the classification of finite-time singularities. Section \ref{sec:sufficient} introduces the angular cone, proves its forward invariance, and establishes a sufficient condition for finite-time collision together with an explicit linear time bound. Section \ref{subsec:quadratic-bound} develops a quadratic refinement and provides its analytical justification for uniform density. Finally, Section \ref{sec:numerical-experiments} presents numerical experiments for different density profiles and initial conditions, including numerical comparisons of the quadratic bound for a uniform segment.


\section{Hamiltonian formalism}
\label{sec:model}

Our formulation follows the notation and methodology developed in \cite{MartinezThesis} for a linear density. In the present work, we extend that framework to an arbitrary positive continuous density, while retaining the same notation whenever possible. 

\begin{figure}[h]
\vspace{-0.60cm}\includegraphics[scale=0.46]{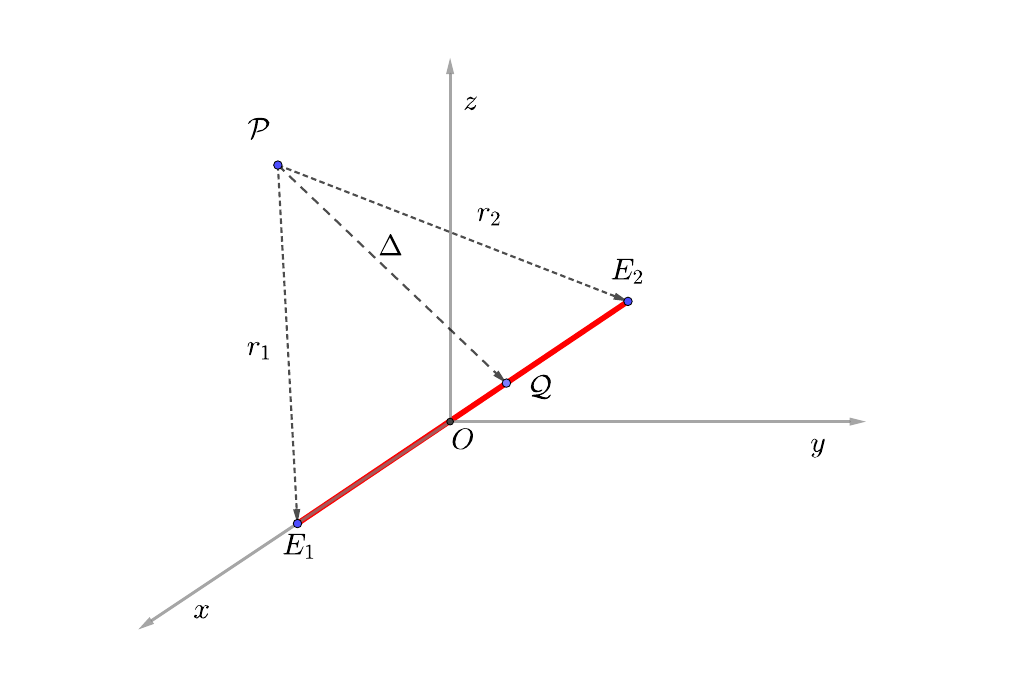}
\hspace{-1.5cm}\includegraphics[scale=0.5]{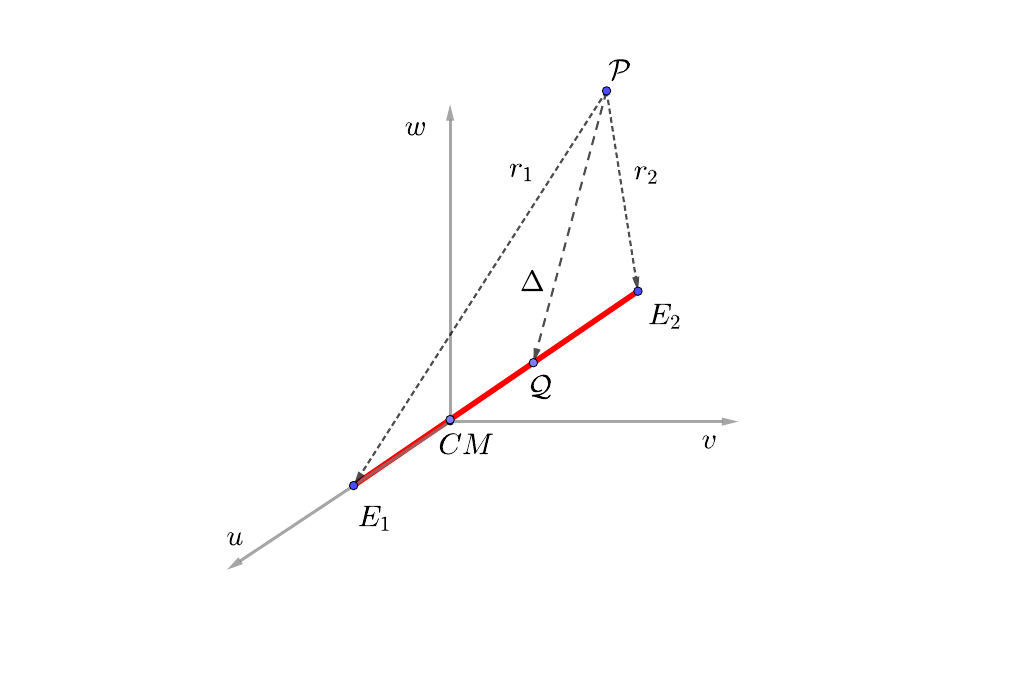}
\vspace{-0.5cm}\caption{\footnotesize Left: The reference frame $\{O;x,y,z\}$ with the segment centered on the origin and lying along the $x$-axis. Right: The new frame $\{O;\xi,\eta,\zeta\}$ {\color{black} resulting from }a translation in the $x$-axis.}
\label{fig:SegmentoFijoPosicionCentrado}
\end{figure}

We fix an inertial reference frame $\{O;x,y,z\}$ with the $2L$-long segment $\mathcal S=\{(x,0,0): x\in [-L,L]\}$ lying along the $x$-axis and centered at the origin, as shown in the left image of Figure~\ref{fig:SegmentoFijoPosicionCentrado}. In this frame, we assume that the density is given by an arbitrary continuous, positive  function 
$$0<\sigma_{min}\leq\sigma(x)\leq\sigma_{max},\quad  x \in [-L, L].$$
We use the following notation to refer to the end points of the segment $E_1 = (L, 0, 0)$ and $E_2 = (-L, 0, 0)$. The vectors joining the massless particle $\mathcal{P} = (x, y, z)$ to the end points are defined as
$$\mathbf{r}_1 = \vec{\mathcal{P} E_1}, \quad \mathbf{r}_2 = \vec{\mathcal{P} E_2}.$$
Here, boldface characters are used to distinguish vectors from their magnitudes, which are expressed in regular characters. Similar to the case of constant density \cite{Riaguas1999}, or linear density \cite{Elisa2024a}, the analysis for arbitrary density will rely on the following distances
$$\Delta = \|\vec{\mathcal{P} \mathcal{Q}}\|, \quad r_1 = \|\vec{\mathcal{P} E_1}\|, \quad r_2 = \|\vec{\mathcal{P} E_2}\|, $$
where $\mathcal{Q} \in \mathcal{S}$ is an arbitrary point on the segment.

The computation of the potential function is done in  a new frame $\{O;\xi,\eta,\zeta\}$ resulting from a translation in the $x$-axis locating the center of mass $(\bar c,0,0)$ at the new origin, see the right image of Figure~\ref{fig:SegmentoFijoPosicionCentrado}. To avoid unnecessary notation, we use the same symbols for the endpoints and the particle in both reference frames, with the understanding that their coordinates are expressed in the new frame. With this agreement, the distances from $\mathcal{P}$ to the end points are
\begin{equation}
 r_{1}= \sqrt{(-L +\xi+\bar{c})^2+\eta^2 + \zeta^2},\quad  r_{2}= \sqrt{(L + \xi+ \bar{c})^2+\eta^2 + \zeta^2},
\end{equation}
and $\mathcal{S}$ is redefined in this new reference frame as $\mathcal{S}=\{(\xi,0,0):\, \xi\in[-L-\bar{c},L-\bar{c}]\}$. Hence, the potential function is given by the following line integral
\begin{equation}
  \label{eq:PotentialFunction}
\begin{array}{rl}
U(\mathcal{P}) & =  -G \displaystyle\int_\mathcal{S} \dfrac{d{\textrm M}}{\Delta}= -G \displaystyle\int_{-L-\bar c}^{L-\bar c} \dfrac{\sigma_c(\chi)}{\|\vec{ \mathcal{P} \mathcal{Q}}\|} \,d\chi, \end{array}
\end{equation}
where $G$ is the universal gravitational constant, and $\sigma_c(\chi):=\sigma(\chi+\bar{c}),\; \chi\in[-L-\bar{c},L-\bar{c}]$ is the density in the translated coordinate.

For the Hamiltonian formulation of the problem, we will employ the following notation for the position and linear momentum of the massless particle $\mathcal{P}$
\[
    \mathbf{q}=(\xi,\eta,\zeta)\in\R^3-\mathcal S,
    \qquad
    \mathbf{p}=(p_\xi,p_\eta,p_\zeta)\in\R^3.
\]
Hence, the Hamiltonian function is given by
\begin{equation}
    \Hcal(\mathbf{q},\mathbf{p})
    =\frac12\left(p_\xi^2+p_\eta^2+p_\zeta^2\right)+U(\mathbf{q}),
    \qquad
    (\mathbf{q},\mathbf{p})\in T^*(\R^3\setminus\Sseg).
    \label{eq:Hamiltonian}
\end{equation}

To simplify the denominator, we express {\color{black} $\vec{ \mathcal{P} \mathcal{Q}}= \mathbf{r}_1 + t\,\mathbf{r}_{12}$, where $\mathbf{r}_{12}:= \vec{E_1 E_2}=\mathbf{r}_2-\mathbf{r}_1$, ${r}_{12}= |\vec{E_1 E_2}|=2L$}  with $t\in(0,1)$, which leads us to $\|\vec{ \mathcal{P} \mathcal{Q}}\|=\sqrt{r_1^2 + t^2\,r_{12}^2+2t\, \mathbf r_1\cdot  \mathbf r_{12}}$.  The above integral is computed, in terms of the new variable $t$, by considering the change of variable $\chi=2Lt-L- \bar{c}$, $d\chi=2Ldt$. Thus, considering the relation {\color{black} $2\, \mathbf r_1\cdot \mathbf r_{12}=r_2^2 -r_1^2 -4L^2$ } obtained from the Cosine Theorem, we are left with the following integral
\begin{equation}
\label{eq:PotentialFunctionBis}
U(\mathbf{q})  = -G \displaystyle\int_0^1 \frac{\hat\sigma(t) }{\sqrt{t^2+\alpha\, t+\beta}} \ dt,
\end{equation}
where the density function is now expressed as $\hat\sigma(t)=\sigma_c\left(\chi(t)\right)$, with $t\in[0,1]$. Moreover, the parameters $\alpha= (-4 L^2-{r}_{1}^2+{r}_{2}^2)/4 L^2$ and $\beta= {r}_{1}^2/4 L^2$. 

\begin{remark}
\label{remark:PotentialFunction}
    The expression of the potential function given in \eqref{eq:PotentialFunctionBis} shows that $U(\mathbf{q}) $  depends exclusively on the distances $r_1$ and $r_2$ of the particle to the endpoints of the segment, and the evaluation of the density function along the segment.
\end{remark}


\begin{remark}
Every trajectory with conserved energy $h<0$ satisfies
\[
    U(\mathbf{q}(t))
    =
    h-\frac12\|\mathbf{p}(t)\|^2
    \leq h<0.
\]
Hence $\mathbf{q}(t)$ remains in the bounded sublevel set $ \{\mathbf{q}:|U(\mathbf{q})|\geq |h|\}$, and therefore cannot escape to spatial infinity. This conclusion does not
exclude collisions with the segment.
\end{remark}


\subsection{Cylindrical Hamiltonian reduction}
\label{subsec:cylindrical}
The simplified expression of the potential function given in \eqref{eq:PotentialFunctionBis}, and Remark~\ref{remark:PotentialFunction} show that $U(\mathbf{q}) $ is invariant under rotations about the $\xi$-axis. Therefore, regardless of the density function that we consider, the Hamiltonian system defined by \eqref{eq:Hamiltonian} admits the following axial-reduction performed through cylindrical coordinates.

On the region $\rho>0$, we define cylindrical coordinates around the segment axis through the symplectic transformation $(\xi,\eta,\zeta,p_\xi,p_\eta,p_\zeta)\longrightarrow (\xi,\rho,\varphi,p_\xi,p_\rho,p_\varphi)$, where 
\begin{equation}
\label{eq:cylindrical}
\begin{aligned}
\eta &= \rho\cos\varphi,\quad &p_\eta &= p_\rho\cos\varphi
-\frac{p_\varphi}{\rho}\sin\varphi,\\
\zeta &= \rho\sin\varphi,\quad &p_\zeta &= p_\rho\sin\varphi
+\frac{p_\varphi}{\rho}\cos\varphi,\\
\rho &= \sqrt{\eta^2+\zeta^2},\quad &p_\varphi &= \eta p_\zeta-\zeta p_\eta
= \rho^2\dot{\varphi}.
\end{aligned}
\end{equation}

As a consequence of axial symmetry, the momentum $p_\varphi$ is a conserved quantity. We will denote its constant value by $\ell:=p_\varphi$. In these coordinates, the Hamiltonian is
\begin{equation}
    \Hcal=\frac12\left(p_\xi^2+p_\rho^2+\frac{\ell^2}{\rho^2}\right)
    +U(\xi,\rho),
    \label{eq:cylindrical-Hamiltonian}
\end{equation}
where $r_{1}= \sqrt{(-L +\xi+\bar{c})^2+ \rho^2},$ $r_{2}= \sqrt{(L + \xi+ \bar{c})^2+\rho^2}$, and Hamilton's equations become
\begin{equation}
\begin{aligned}
    \dot\xi&=p_\xi,
    &\dot p_\xi&=-U_\xi(\xi,\rho),\\
    \dot\rho&=p_\rho,
    &\dot p_\rho&=\frac{\ell^2}{\rho^3}-U_\rho(\xi,\rho),\\
    \dot\varphi&=\frac{\ell}{\rho^2},
    &\dot\ell&=0.
\end{aligned}
\label{eq:cylindrical-equations}
\end{equation}

\section{Finite-time singularities}
\label{sec:singularities}

This section shows that every finite-time singularity forces the trajectory
to approach the segment regardless of the density function we assume. Before stating this result, we recall some basic terminology and notation relative to singularities.

Let $\mathbf{z}_0=(\mathbf{q}_0,\mathbf{p}_0)\in T^*(\R^3\setminus\Sseg)$ be an initial condition. We denote corresponding solution of the Hamiltonian system defined by \eqref{eq:Hamiltonian} as $\mathbf z(t)=\Phi(t,\mathbf{z}_0)$. This solution is defined on a maximal interval $(\omega_-,\omega_+)$. A finite endpoint of this interval is called a singular time, and the solution $\mathbf z(t)$ is said to be singular.

\begin{definition}
A singularity at $\omega_+$ is called a collision singularity if there exists $\mathbf{q}_*\in\Sseg$ such that
\begin{equation}
    \lim_{t\to\omega_+}\mathbf{q}(t)=\mathbf{q}_*.
    \label{eq:collision-definition}
\end{equation}
The corresponding definition at $\omega_-$ is analogous with the time orientation reversed.  In this case, the singularity is called ejection singularity.
\end{definition}

The potential is smooth on $\R^3\setminus\Sseg$ but is unbounded below at the segment. The next lemma records the singular behavior needed below providing bounds of the potential.

\begin{lemma}
\label{Lem:PotentialBounds}
For an arbitrary positive and bounded density distribution along the
segment, the potential function \eqref{eq:PotentialFunctionBis} satisfies
\begin{equation}
    k\log\left(
    1+\frac{2L}{\delta(\mathbf{q})}
    \right)\leq|U(\mathbf{q})|
    \leq
    2K\log\left(
    1+\frac{2L}{\delta(\mathbf{q})}
    \right),
    \label{eq:PotentialBounds}
\end{equation}
where $K=G\,\sigma_{\max},$ $k=G\,\sigma_{\min},$ and $\delta(\mathbf{q})=\operatorname{dist}(\mathbf{q},\mathcal{S})$ is the distance from the particle to the segment.
\end{lemma}

\begin{proof}
From the definition of the density function $\hat\sigma$ in
\eqref{eq:PotentialFunctionBis}, which is a reparameterization of the
density $\sigma$ along the segment, we have $0< \sigma_{min}\leq \hat\sigma(t)\leq \sigma_{\max}$, with $t\in[0,1]$.
Therefore,
\begin{align}
    k
    \log\left(
    \frac{s+2L}{s-2L}
    \right)=\int_0^1
    \frac{G\,\sigma_{\min}}
    {\sqrt{t^2+\alpha t+\beta}}
    \,\dd t\leq|U(\mathbf{q})|
    &\leq
    \int_0^1
    \frac{G\,\sigma_{\max}}
    {\sqrt{t^2+\alpha t+\beta}}
    \,\dd t =K
    \log\left(
    \frac{s+2L}{s-2L}
    \right),
    \label{eq:first-potential-bound}
\end{align}
where $s=r_1+r_2$. Lower and upper bounds are  the well known expression for the potential of a segment with constant density. Now, we will relate these bounds with $\delta(\mathbf{q})$, i. e., the distance from $\mathbf{q}$ to the segment. This part will be done by proving two inequalities.

\vspace{0.4cm}
\noindent\textit{First inequality:} Let $\mathbf{q}^{*}\in\mathcal{S}$ be a point realizing the distance to the
segment, so that $ \|\mathbf{q}-\mathbf{q}^{*}\|=\delta(\mathbf{q}).$ By the triangle inequality,
\[
\begin{aligned}
    s=r_1+r_2
    &\leq
    \|\mathbf{q}-\mathbf{q}^{*}\|
    +\|\mathbf{q}^{*}-E_1\|
    +\|\mathbf{q}-\mathbf{q}^{*}\|
    +\|\mathbf{q}^{*}-E_2\| =2\delta(\mathbf{q})+2L.
\end{aligned}
\]
Consequently,
\begin{equation}
     s +2L
    \leq 2\delta(\mathbf{q})+4L.
    \label{eq:sum-upper-bound}
\end{equation}

\vspace{0.4cm}
\noindent\textit{Second inequality:} We next prove that
\begin{equation}
     s ^2-4L^2
    \geq 4\delta(\mathbf{q})^2.
    \label{eq:sum-lower-bound}
\end{equation}
For this part, we translate the reference frame and assume that the segment has endpoints
$(-L,0,0)$ and $(L,0,0)$. Moreover, we restrict to the plane containing the segment and the particle, and employ the notation $\rho:=\sqrt{\eta^2+\zeta^2}$ to denote the normal component of the particle to the segment.

Suppose first that the orthogonal projection of $\mathbf{q}$ onto the
$\xi$-axis lies inside the segment, that is, $|\xi|\leq L$. In this case,
$\delta(\mathbf{q})=\rho$, $r_1=\sqrt{(L-\xi)^2+\rho^2}$, and $r_2=\sqrt{(L+\xi)^2+\rho^2}$. Applying the triangle inequality to the vectors $\mathbf{r}_1=(L-\xi,\rho)$ and $\mathbf{r}_2=(L+\xi,\rho)$, we obtain
\[
\begin{aligned}
    s
    &\geq
    \sqrt{(2L)^2+(2\rho)^2} =2\sqrt{L^2+\delta(\mathbf{q})^2}.
\end{aligned}
\]
Squaring the last inequality yields \eqref{eq:sum-lower-bound}.

If the projection lies outside the segment, we distinguish two cases: $\xi\geq L$ and $\xi\leq - L$. We will prove the first case, since the second is analogous. Hence, assuming that $\xi\geq L$, let us define $d=\xi-L\geq0$. Then, the distance is given as
\[
    \delta(\mathbf{q})=\sqrt{d^2+\rho^2}.
\]
Applying again the triangle inequality to the vectors $\mathbf{r}_1=(d,\rho)$ and $\mathbf{r}_2=(d+2L,\rho)$ gives
\[
    s
    \geq
    2\sqrt{(d+L)^2+\rho^2}
    \geq
    2\sqrt{L^2+\delta(\mathbf{q})^2}.
\]
Thus, squaring the last inequality yields \eqref{eq:sum-lower-bound}, which
holds in all cases.

Finally, we prove \eqref{eq:PotentialBounds}  using the inequalities \eqref{eq:sum-upper-bound} and \eqref{eq:sum-lower-bound}. Precisely, we consider the following bounds for the argument of the logarithm in \eqref{eq:first-potential-bound}
\begin{gather}
\begin{aligned}
    \frac{s+2L}{s-2L}
    =&
    \frac{\bigl(s+2L\bigr)^2}
    {s^2-4L^2} \leq
    \frac{\bigl(2\delta(\mathbf{q})+4L\bigr)^2}
    {4\delta(\mathbf{q})^2}= \left(
    1+\frac{2L}{\delta(\mathbf{q})}
    \right)^2.\\
     \frac{s+2L}{s-2L}
    =&
    1+\frac{4L}
    {s -2L} \geq
    1+\frac{2L}
    {\delta(\mathbf{q})}.
    \label{eq:ratio-distance-bound}
\end{aligned}
\end{gather}
Then, the proof is completed by taking logarithms in \eqref{eq:ratio-distance-bound}.
\end{proof}


The singularity of the potential is the only obstruction to applying the standard continuation theorem to the Hamiltonian vector field. Moreover, Lemma~\ref{Lem:PotentialBounds} showed that any finite-time singularity requires the distance $\delta(\mathbf{q})$ from the particle to the segment to tend to zero. These
observations lead to the following theorem.

\begin{theorem}[Finite-time singularities]
\label{thm:finite-time-singularities}
Let $\mathbf z(t)=(\mathbf{q}(t),\mathbf{p}(t))$ be a solution of the
Hamiltonian system defined by \eqref{eq:Hamiltonian} with maximal interval
$I_{\mathbf{z}_0}=(\omega_-,\omega_+)$. If $\omega_+<+\infty$, then
\begin{equation}
    \lim_{t\to\omega_+}\delta(\mathbf{q}(t))=0.
    \label{eq:finite-singularity-approach}
\end{equation}
The analogous statement holds if $\omega_->-\infty$.
\end{theorem}

\begin{proof}

We will prove the statement for a finite forward endpoint $\omega_+<+\infty$. The finite backward endpoint case is analogous.

Let us assume that $(t_0,\mathbf{z_0})$ is an initial condition determining the solution $\mathbf{z}(t)$ in the energy level $h_0:=\mathcal H(\mathbf z_0)$. Moreover, we assume that  $\omega_+<+\infty$, but $\delta(\mathbf q)$ does not converges to zero as $t\rightarrow\omega_+$. Thus, there exists $\varepsilon_0>0$, a tubular
neighborhood
\[
    \mathcal{N}_{\varepsilon_0}
    :=
    \left\{
        \mathbf q\in\mathbb{R}^3:
        \delta(\mathbf q)\leq\varepsilon_0
    \right\},
\]
and a sequence $t_n\rightarrow\omega_+$ such that $\mathbf{q}_n:=\mathbf{q}(t_n)\notin \mathcal{N}_{\varepsilon_0}$ for all $n\in\mathbb{N}$.

Since $\Sseg$ is compact, $\mathcal{N}_{\varepsilon_0}$ is bounded, and hence compact. Moreover,
Lemma~\ref{Lem:PotentialBounds} implies that, whenever
$\mathbf q\notin\mathcal{N}_{\varepsilon_0}$,
\[
    U(\mathbf q)
    \geq
    -2K\log\left(
        1+\frac{2L}{\varepsilon_0}
    \right)
    =:-C_{\varepsilon_0}.
\]
By energy conservation we can ensure bounded velocity. Namely, using the previous inequality, we have 
$$\|\mathbf p(t)\|^2=2h_0-2U(\mathbf q(t))\leq 2h_0+2C_{\varepsilon_0},$$ 
whenever $\mathbf q\notin\mathcal{N}_{\varepsilon_0}$. Moreover, since $\dot{\mathbf q}=\mathbf p$ is bounded, and $\omega_+<+\infty$, the total length of the curve arcs outside of $\mathcal{N}_{\varepsilon_0}$ is finite. 
Thus, the portion of $\mathbf{z}(t)$ with corresponding $\mathbf{q}(t)$ outside $\mathcal{N}_{\varepsilon_0}$ lies in a compact subset $K\subset T^*(\mathbb{R}^3\setminus\Sseg)$. 

The Hamiltonian vector field is smooth on
$T^*(\mathbb{R}^3\setminus\Sseg)$, and in the compact set $K$. Therefore, by the local existence and
continuation theorem, there exists $\tau>0$, independent of $n$, such that
the solution with initial condition $\mathbf z(t_n)$ exists on
$[t_n,t_n+\tau]$. For $n$ sufficiently large,
$t_n+\tau>\omega_+$. By uniqueness, this solution agrees with the original
one on $[t_n,\omega_+)$ and consequently extends it beyond $\omega_+$,
contradicting the maximality of $\omega_+$.
\end{proof}


This theorem says that the model has no finite-time breakdown mechanism other than approach to the material segment. It  was first proven for linear density in \cite{MartinezThesis}, and extended here to arbitrary density. However, the more informative question is which conditions actually produce such an approach.

\section{Necessary conditions for collision}
\label{sec:CollisionConditions}

\begin{theorem}[Necessary angular-momentum condition]
\label{thm:necessary-angular-momentum}
If a solution satisfies
\[
\lim_{t\rightarrow\omega_+}\delta(\mathbf q(t))=0,
\]
then its axial angular momentum \(\ell\) vanishes. Therefore, every solution that approaches the segment $\Sseg$ as $t\rightarrow\omega_+$ is confined to a fixed plane containing the segment.
\end{theorem}

\begin{proof}
Let $h$ be the conserved value of the Hamiltonian. Suppose, by
contradiction, that $\ell\neq0$. Along any sequence
$\mathbf{q}_n\to\Sseg$, approach to the segment implies $\rho_n\longrightarrow0$ with $ \rho_n=\sqrt{\eta_n^2+\zeta_n^2}$. By Lemma~\ref{Lem:PotentialBounds}, the potential diverges at most
logarithmically in the distance to the segment $\delta(\mathbf q_n)$, whereas
$\ell^2/(2\rho_n^2)$ diverges quadratically in $1/\rho_n$. Hence, considering $\delta(\mathbf q_n)\geq \rho_n$
\begin{equation}
    \mathcal{E}_n=\frac{\ell^2}{2\rho_n^2}+U(\mathbf{q}_n)
    \longrightarrow+\infty,
    \label{eq:centrifugal-divergence}
\end{equation}
and the evaluation of the Hamiltonian \eqref{eq:cylindrical-Hamiltonian} along the trajectory approaching the segment is given by $\mathcal H(\mathbf{q}_n,\mathbf{p}_n)=1/2(p_\xi^2+p_\rho^2)+ \mathcal{E}_n$, which also diverges to $+\infty$. This contradicts the finiteness and conservation of the energy.

Moreover, $\ell=0$ implies $\varphi(t)=\varphi_0$ is constant, and the trajectory remains in the plane $\Pi_{\varphi_0}:=\{\varphi=\varphi_0\}$, which contains the $\xi$-axis and, consequently, the segment $\Sseg$.
\end{proof}


By Theorem~\ref{thm:necessary-angular-momentum}, every collision trajectory is contained in a
fixed plane containing the segment. Since the Hamiltonian is invariant under
rotations about the $\xi$-axis, all such planes are dynamically equivalent.
We therefore fix an arbitrary invariant plane $\Pi_{\varphi_0}$ and identify
its upper half-plane $\Pi_{\varphi_0}^+$ with coordinates $(\xi,\rho)$, where $\rho>0$.

Moreover, in what remains we will employ the following notation. In the translated reference frame, let the endpoints of the segment be given by $\xi_-:=-L-\bar c$, $\xi_+:=L-\bar c$, and denote the material segment in the invariant plane by $\mathcal{S}=I\times\{0\}$, where $I=[\xi_-,\xi_+].$ For $\chi\in I$, let $\mathbf{s}_{\chi}:=(\chi,0)\in\mathcal{S}.$

Additionally, with a slight abuse of notation we will employ the same symbols $(\mathbf{q},\mathbf{p})$ for the axial-symmetry reduced coordinates in the invariant half-plane $\Pi_{\varphi_0}^+$. That is to say, from now on, we have
\begin{equation}
    \label{eq:PlanarCoordinates}
\mathbf{q}=(\xi,\rho),\qquad \mathbf{p}=(p_\xi,p_\rho).
\end{equation}

\begin{remark}
    Note that Theorem~\ref{thm:necessary-angular-momentum} does not exclude the possibility of colliding to the segment through a straight line contained in the equatorial plane when the density distribution is symmetric with respect to the midpoint of the segment. These solutions are analogous to the collisions detected in \cite{AlbertiVidal2010,AlbertiVidal2019}.
\end{remark}

\begin{remark}
    \label{rem:Acceleration}
     The acceleration $\ddot{\mathbf{q}}(t)$ of any particle subjected to the gravitational pull of the segment is a vector pointing from
$\mathbf{q}(t)$ toward the material segment. Precisely, since $\ell=0$ on an invariant plane, the Hamiltonian restricted to
$\Pi_{\varphi_0}$ is
\begin{equation}
    \Hcal_2(\xi,\rho,p_\xi,p_\rho)
    =
    \frac12\left(p_\xi^2+p_\rho^2\right)
    +U(\xi,\rho),
    \qquad \rho>0,
    \label{eq:invariant-plane-Hamiltonian}
\end{equation}
and the equations of motion are
\begin{equation}
    \dot{\mathbf{q}}=\mathbf{p},
    \qquad
    \dot{\mathbf{p}}=-\nabla U(\mathbf{q}).
    \label{eq:invariant-plane-equations}
\end{equation}
Differentiating the line-integral representation of the potential gives
\begin{equation}
\begin{aligned}
    \ddot{\mathbf{q}}
    =\dot{\mathbf{p}}
    =-\nabla U(\mathbf{q})
    &=
    G\int_{I}
    \frac{{\sigma_c}(\chi)}
    {\left((\xi-\chi)^2+\rho^2\right)^{3/2}}
    \bigl(\mathbf{s}_{\chi}-\mathbf{q}\bigr)\,\dd\chi .
\end{aligned}
\label{eq:planar-gravitational-acceleration}
\end{equation}
Thus, the acceleration is a positive combination of vectors pointing from
$\mathbf{q}$ toward the interior of the material segment. Moreover, only when $\rho=0$ we have the acceleration as a vector along the line containing the segment.
\end{remark}


\begin{theorem}[Approach to the segment implies collision]
\label{theorem:approach-implies-collision}
Let $\mathbf z(t)=(\mathbf q(t),\mathbf p(t))$ be a planar solution ($\ell=0$) with maximal forward interval
$I_{\mathbf z_0}^+=[0,\omega_+)$, where $\omega_+<+\infty$. Assume that
\[
\lim_{t\rightarrow\omega_+}\delta(\mathbf q(t))=0.
\]
Then there exists a point $\mathbf q_*\in\mathcal{S}$ such that
\[
\lim_{t\rightarrow\omega_+}\mathbf q(t)=\mathbf q_*.
\]
Consequently, the singularity at $\omega_+$ is a collision singularity.
\end{theorem}

\begin{proof}
For $\chi\in I$, we define
\begin{equation}
    \label{eq:Kappa}
    w_{\mathbf q}(\chi)
={G{\sigma_c}(\chi)}/
{\big((\xi-\chi)^2+\rho^2\big)^{3/2}},\qquad
   \kappa(\mathbf q)=\int_I w_{\mathbf q}(\chi)\,d\chi.
\end{equation}
Since $\delta(\mathbf q(t))\geq \rho(t),$ the hypothesis implies that $\rho(t)\to0$ as
$t\rightarrow\omega_+$. Additionally, in the invariant plane, we have $\dot p_\rho=-\kappa(\mathbf q)\rho<0$. Since $\rho(t)\to0$ while $\rho(t)>0$ for $t<\omega_+$, there exists
$t_1\in(0,\omega_+)$ such that $p_\rho(t_1)<0$. Moreover, since $p_\rho$ is
strictly decreasing, there exists $c>0$ such that
\[
p_\rho(t)\leq -c<0,
\qquad
t\in[t_1,\omega_+).
\]
Consequently, $\rho$ is strictly decreasing on this interval and
\[
-\dot\rho(t)\geq c.
\]

Let $h:=\mathcal H_2(\mathbf{z}_0)$ denote the conserved energy. 
Lemma~\ref{Lem:PotentialBounds} and $\delta(\mathbf q)\geq\rho$ imply
\[
|p_\xi(t)|
\leq
\sqrt{
2|h|+
4K\log\left(1+\frac{2L}{\rho(t)}\right)
}.
\]
Using $-\dot\rho(t)\geq c$, we obtain
\[
\begin{aligned}
\int_{t_1}^{\omega_+}|p_\xi(t)|\,\dd t
&\leq
\frac{1}{c}
\int_0^{\rho(t_1)}
\sqrt{
2|h|+
4K\log\left(1+\frac{2L}{r}\right)
}
\,\dd r
<+\infty.
\end{aligned}
\]
Moreover,
\[
\int_{t_1}^{\omega_+}|p_\rho(t)|\,\dd t
=
-\int_{t_1}^{\omega_+}p_\rho(t)\,\dd t
=
\rho(t_1)<+\infty.
\]
Since the solution is regular on the compact interval $[0,t_1]$, it follows
that
\[
\int_0^{\omega_+}\|\dot{\mathbf q}(t)\|\,\dd t<+\infty.
\]
Thus, the trajectory has finite length. In particular, for
$s<t<\omega_+$,
\[
\|\mathbf q(t)-\mathbf q(s)\|
\leq
\int_s^t\|\dot{\mathbf q}(u)\|\,\dd u,
\]
and the right-hand side tends to zero as $s,t\rightarrow\omega_+$. Hence, there exists $\mathbf q_*\in\mathbb R^3$ such that
\[
\lim_{t\rightarrow\omega_+}\mathbf q(t)=\mathbf q_*.
\]
Finally, since $\delta(\mathbf q(t))\to0$ and $\mathcal S$ is closed,
we have $\delta(\mathbf q_*)=0$, and therefore
$\mathbf q_*\in\mathcal S$.
\end{proof}


\begin{corollary}
    Every finite-time singularity of the system is either a collision singularity or an ejection singularity.
\end{corollary}
\begin{proof}
The forward-time statement follows from
Theorem~\ref{thm:finite-time-singularities}
    and Theorem~\ref{theorem:approach-implies-collision}. The backward-time statement follows by applying the same argument to the time-reversed solution.
\end{proof}

\section{Sufficient conditions for collision}
\label{sec:sufficient}

\begin{definition}
\label{def:AngularCone}
For every point $\mathbf{q}_*\in\Pi_{\varphi_0}^+$ in configuration space, we define the $\mathbf{q}_*$-triangle as
\[
    \mathcal{T}(\mathbf{q}_*)
    :=
   \Bigl\{ 
            \mathbf{x} \in \Pi_{\varphi_0}^+:\,  \mathbf{x} = \mathbf{q}_* + \lambda_1\bigl(\mathbf{s}_{\chi_1}-\mathbf{q}_*\bigr), \, \chi_1\in I,\ \lambda_1 \in [0,1) 
        \Bigr\}
\]
Moreover, for every $\mathbf{x}\in\mathcal{T}(\mathbf{q}_*)$, the set $ {PWC}(\mathbf{x})
    :=
    \left\{\mathbf{y}\in\mathbb{R}^2:\,\mathbf{y}=
        \lambda_2\bigl(\mathbf{s}_{\chi_2}-\mathbf{x}\bigr),\,
        \chi_2\in I,\ \lambda_2\geq0
    \right\}$ is the \emph{pointwise cone based at} $\mathbf{x}$. Then, we define the \emph{angular cone associated with} $\mathbf{q}_*$ as the set
\begin{equation}
    \mathcal{C}(\mathbf{q}_*)
    :=
    \left\{
        (\mathbf{x},\mathbf{y}):
        \mathbf{x}\in\mathcal{T}(\mathbf{q}_*),
        \ \mathbf{y}\in{PWC}(\mathbf{x})
    \right\}.
    \label{eq:angular-cone}
\end{equation}
Equivalently, $(\mathbf{x},\mathbf{y})\in\mathcal{C}(\mathbf{q}_*)$ if and only if there exist $\chi_1,\chi_2\in I$ and
$\lambda_1\in[0,1)$, $\lambda_2\geq0$ such that
\[
    \mathbf{x}
    =
    \mathbf{q}_*
    +\lambda_1\bigl(\mathbf{s}_{\chi_1}-\mathbf{q}_*\bigr),
    \qquad
    \mathbf{y}
    =
    \lambda_2\bigl(\mathbf{s}_{\chi_2}-\mathbf{x}\bigr),
    \qquad
    \mathbf{x}\notin\mathcal{S}.
\]
\end{definition}

The set $\mathcal{C}(\mathbf{q}_*)$ encodes both the configuration region and the corresponding directions over each point. Its projection onto the configuration space is $ \mathcal{T}(\mathbf{q}_*)$. For each $\mathbf{x}\in\mathcal{T}(\mathbf{q}_*)$, the associated fiber is the pointwise cone $PWC(\mathbf{x})$. When this construction is applied to a trajectory, the second component is interpreted according to the dynamical quantity under consideration. Using the Euclidean identification $\mathbb{R}^2\simeq(\mathbb{R}^2)^*$, the same admissibility condition can therefore be written for different nature pairs as
\[
    (\mathbf{q},\dot{\mathbf{q}})\in\mathcal{C}(\mathbf{q}_*),
    \qquad
    (\mathbf{q},\mathbf{p})\in\mathcal{C}(\mathbf{q}_*),
    \qquad
    (\mathbf{q},\dot{\mathbf{p}})\in\mathcal{C}(\mathbf{q}_*),
     \qquad\text{or}\qquad
     (\mathbf{q},\ddot{\mathbf{q}})\in\mathcal{C}(\mathbf{q}_*),
\]
regardless of whether the relevant vector is the velocity, momentum, or
momentum derivative.

{\color{black}

\begin{lemma}[Forward invariance of the angular cone]
\label{lem:tecnico1}
Let $\mathbf z_0=(\mathbf q_0,\mathbf p_0)=(\xi_0,\rho_0,p_{\xi_0},p_{\rho_0})$ with $p_{\rho_0}<0$ an initial condition in the invariant upper half-plane $\Pi_{\varphi_0}^+$ satisfying $(\mathbf q_0,\mathbf p_0)\in \mathcal C(\mathbf q_0)$, and let $(\mathbf q(t),\mathbf p(t))$ be the corresponding solution of system \eqref{eq:invariant-plane-equations} with  maximal forward interval $I_{\mathbf{z}_0}^+=[0,\omega_+)$. Then $(\mathbf q(t),\mathbf p(t))\in \mathcal C(\mathbf q_0)$ for all $t\in I_{\mathbf{z}_0}^+$.
\end{lemma}

\begin{proof} For the sake of compact expressions, we eliminate the explicit reference to the dependence of the independent variable $t$ when unnecessary. Hence we write $(\mathbf q,\mathbf p)$ instead of $(\mathbf q(t),\mathbf p(t))$ and we apply the same criteria for every $t$-dependent function.

To show that $(\mathbf q(t),\mathbf p(t))\in \mathcal C(\mathbf q_0)$, we need to prove that $\mathbf q(t)\in \mathcal T(\mathbf q_0)$ and $\mathbf p(t)\in PWC(\mathbf q(t))$ for all $t\in I_{\mathbf{z}_0}^+$. 

\noindent \textit{Part I, $\mathbf p\in PWC(\mathbf q)$:} Using $w_{\mathbf q}$ and $ \kappa(\mathbf q)$ given in \eqref{eq:Kappa}, we define
\[
\overline{\chi}(\mathbf q)
=
\frac{\int_I \chi w_{\mathbf q}(\chi)\,d\chi}
{\int_I w_{\mathbf q}(\chi)\,d\chi}.
\]
Since ${\sigma}_c>0$, we have $\kappa(\mathbf q)>0$ and $\overline{\chi}(\mathbf q)\in(\xi_-,\xi_+)$. Hence, using this notation the gradient of the potential function $-\nabla U=[\dot p_\xi,\dot p_\rho]$ becomes
\begin{equation}
    \label{eq:GradPotencial}
    \dot p_\xi
=
\kappa(\mathbf q)\bigl(\overline{\chi}(\mathbf q)-\xi\bigr),
\qquad
\dot p_\rho=-\kappa(\mathbf q)\rho.
\end{equation}
Moreover, \(\dot p_\rho=-\kappa(\mathbf q)\rho<0\) whenever
\(\rho>0\), and \(p_\rho(0)<0\), we have
\(p_\rho(t)<0\) throughout the interval of existence before the first
intersection with $\rho=0$. Hence, define $\lambda:=-{p_\rho}/{\rho}>0$ and $\chi_{\mathbf p}:=\xi+{p_\xi}/{\lambda}$. Then
\[
\mathbf p=\lambda(s_{\chi_{\mathbf p}}-\mathbf q).
\]
Consequently, $\mathbf p\in PWC(\mathbf q)$ iff $\chi_{\mathbf p}\in I$. A direct calculation of the derivative of $\chi_{\mathbf p}$ gives the following initial value problem determining the evolution of $\chi_{\mathbf p}$
\[
\dot\chi_{\mathbf p}
=
\frac{\kappa}{\lambda}
\bigl(\overline{\chi}(\mathbf q)-\chi_{\mathbf p}\bigr),\qquad \chi_{\mathbf p}(0)\in I.
\]
This is a scalar differential equation that attracts $\chi_{\mathbf p}$ toward $\overline{\chi}(\mathbf q)\in I$. Therefore the interval $I$ is forward invariant, $\chi_{\mathbf p}(t)\in I$ for all $t\in I_{\mathbf{z}0}^+$, and hence
\[
\mathbf p(t)\in PWC(\mathbf q(t)),\quad t\in I_{\mathbf{z}0}^+.
\]

\par\noindent
\begin{minipage}[t]{\linewidth}
\begin{wrapfigure}[14]{r}{0.48\linewidth}
    \centering
    \vspace{-1.0cm}
    \begin{tikzpicture}[
        x=1cm,
        y=1cm,
        every node/.style={font=\small},
        boundary/.style={thick},
        projection/.style={densely dashed, gray!70},
        section/.style={densely dashed, thick},
        point/.style={circle, fill=black, inner sep=1.5pt}
    ]

        \coordinate (xmin) at (-2.8,0);
        \coordinate (xmax) at ( 2.8,0);
        \coordinate (qzero) at (0.75,3.6);

        \coordinate (Lrho) at ($(xmin)!0.43!(qzero)$);
        \coordinate (Rrho) at ($(xmax)!0.43!(qzero)$);
        \coordinate (qt)   at (0.45,1.55);

        \coordinate (Lbase) at (Lrho |- xmin);
        \coordinate (Rbase) at (Rrho |- xmax);

        \fill[gray!10] (xmin) -- (xmax) -- (qzero) -- cycle;

        \draw[boundary] (xmin) -- (qzero) -- (xmax);
        \draw[line width=1.2pt] (xmin) -- (xmax);

        \draw[section] (Lrho) -- (Rrho);

        \draw[projection] (Lrho) -- (Lbase);
        \draw[projection] (Rrho) -- (Rbase);

        \node[point] at (xmin) {};
        \node[point] at (xmax) {};
        \node[point] at (qzero) {};
        \node[point] at (Lrho) {};
        \node[point] at (Rrho) {};
        \node[point] at (qt) {};

        \node[below left=2pt]  at (xmin) {$\xi_-$};
        \node[below right=2pt] at (xmax) {$\xi_+$};

        \node[above right=2pt] at (qzero)
            {$q_0=(\xi_0,\rho_0)$};

        \node[left=3pt] at (Lrho)
            {$\bigl(L(\rho),\rho\bigr)$};

        \node[right=3pt] at (Rrho)
            {$\bigl(R(\rho),\rho\bigr)$};

        \node[above=3pt] at (qt)
            {$q(t)=(\xi(t),\rho(t))$};

        \node[below=4pt] at (Lbase) {$L(\rho)$};
        \node[below=4pt] at (Rbase) {$R(\rho)$};


    \end{tikzpicture}
    \vspace{-0.5cm}
    \caption{Geometric construction of functions \(L(\rho)\) and \(R(\rho)\).}
    \label{fig:L-R-region}
\end{wrapfigure}
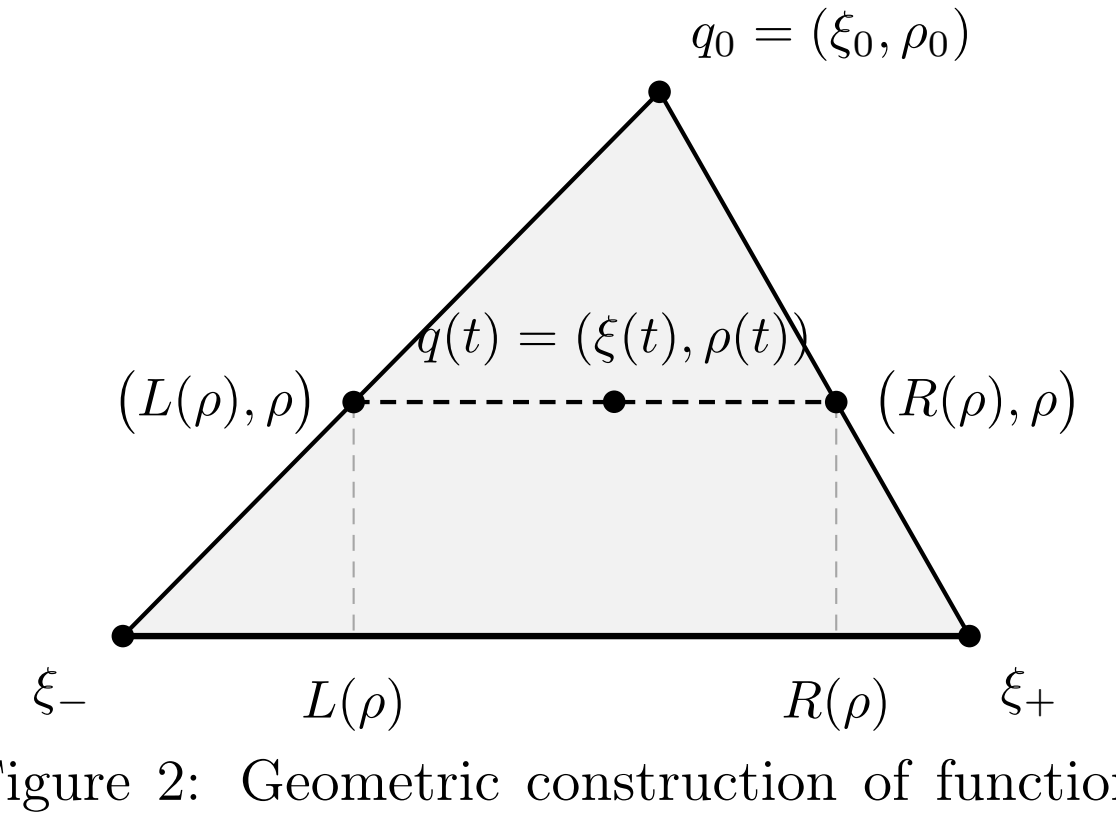

\noindent \textit{Part II, $\mathbf q(t)=(\xi(t),\rho(t))\in \mathcal T(\mathbf q_0)$:} This condition is equivalent to check that $0<\rho<\rho_0$ and $\xi$ is between the bounds determined by the lines joining $\mathbf q_0$ and the endpoints of the segment. Since \(p_\rho(t)<0\) before the first intersection with
\(\rho=0\), and \(\dot\rho=p_\rho\), we have
\(0<\rho(t)\leq\rho_0\). For the second condition we define the bounds for $\xi$ given in Figure~\ref{fig:L-R-region} as 
\[
L(\rho)
=
\alpha\xi_0+(1-\alpha)\xi_-,
\quad
R(\rho)
=
\alpha\xi_0+(1-\alpha)\xi_+.
\]
where $\alpha={\rho}/{\rho_0}$. Then, we define $\mathcal{T}(\mathbf{q}_0)$ as 
\end{minipage}
\par
\[
\mathcal T(\mathbf q_0)
=
\{(\xi,\rho):0<\rho\le\rho_0,\,
L(\rho)\le\xi\le R(\rho)\}.
\]
Thus, $\mathbf q\in \mathcal T(\mathbf q_0)$ is equivalent to prove that the following functions are always positive
\[
F_-=\xi-L(\rho),
\qquad
F_+=R(\rho)-\xi.
\]
To prove that property we consider the initial values $F_-(0)=F_+(0)=0$, and we will prove that \(F_-(t)\geq0\) and \(F_+(t)\geq0\). Precisely, using the definition of $F_-$, $L(\rho)$, and $\alpha$ we obtain
\begin{align*}
\dot F_-
&=\dot{\xi}-L'(\rho)\dot{\rho}=p_\xi-\frac{\xi_0-\xi_-}{\rho_0}\,p_\rho=\lambda(\chi_{\mathbf p}-\xi)
  +\lambda\frac{\rho}{\rho_0}(\xi_0-\xi_-)\\
&=\lambda\left[
\chi_{\mathbf p}-\xi
+\alpha(\xi_0-\xi_-)
\right]=\lambda\left[
\chi_{\mathbf p}-\xi_--\xi
+\alpha(\xi_0-\xi_-)+\xi_-
\right]\\
&=\lambda\left[
\chi_{\mathbf p}-\xi_- -\left(\xi
-(\alpha\xi_0+(1-\alpha)\xi_-)\right)
\right]=\lambda\left[
\chi_{\mathbf p}-\xi_- -(\xi -L(\rho))
\right]
\end{align*}
Proceeding in the same way with $F_+$ leads to the following differential relations
\[
\dot F_-
=
\lambda\bigl(\chi_{\mathbf p}-\xi_- -F_-\bigr),
\qquad
\dot F_+
=
\lambda\bigl(\xi_+-\chi_{\mathbf p}-F_+\bigr),
\]
variation of constants gives
\[
F_-(t)
=
\int_0^t
\lambda(s)\bigl(\chi_{\mathbf p}(s)-\xi_-\bigr)
\exp\left(-\int_s^t\lambda(u)\,du\right)\,ds
\ge0,
\]
and
\[
F_+(t)
=
\int_0^t
\lambda(s)\bigl(\xi_+-\chi_{\mathbf p}(s)\bigr)
\exp\left(-\int_s^t\lambda(u)\,du\right)\,ds
\ge0.
\]
Indeed, for \(0\leq s\leq t\) before the first intersection with
\(\rho=0\), we have \(\lambda(s)>0\) and
\(\chi_{\mathbf p}(s)\in I=[\xi_-,\xi_+]\). Hence
\(\chi_{\mathbf p}(s)-\xi_-\ge0\) and \(\xi_+-\chi_{\mathbf p}(s)\ge0\), while the exponential factor is
strictly positive. Therefore the integrands are nonnegative. Therefore $\mathbf q\in \mathcal T(\mathbf q_0)$.
\end{proof}

}

\begin{theorem}[Sufficient condition for collision]
\label{theorem:SufficientCondition}
Let $\mathbf{z}_0=(\mathbf{q}_0,\mathbf{p}_0)$ be an initial condition in
the invariant upper half-plane, with
$\mathbf{q}_0=(\xi_0,\rho_0)$ and $\rho_0>0$. Let
$\mathbf z(t)=(\mathbf{q}(t),\mathbf{p}(t))$ be the corresponding solution of
\eqref{eq:invariant-plane-equations}. Assume that $(\mathbf{q}_0,\mathbf{p}_0)\in\mathcal{C}(\mathbf{q}_0)$ and $p_{\rho_0}:=p_\rho(0)<0$. Then the forward maximal interval $I_{\mathbf{z}_0}^+=[0,\omega_+)$ associated to $\mathbf{z}(t)$ is bounded, and there exists a point $\mathbf q_*\in\mathcal{S}$ such that
\[
\lim_{t\rightarrow\omega_+}\mathbf q(t)=\mathbf q_*.
\]
Moreover, $\omega_+$ satisfies
\[
\omega_+\leq -\frac{\rho_0}{p_{\rho_0}}.
\]
\end{theorem}
\begin{proof}
By Lemma~\ref{lem:tecnico1}, $\bigl(\mathbf q(t),\mathbf p(t)\bigr)
\in\mathcal C(\mathbf q_0)$ for $t\in I_{\mathbf z_0}^+$. In particular, $\mathbf q(t)$ remains in the invariant upper half-plane,
so $\rho(t)>0$ for $t\in I_{\mathbf z_0}^+$.

By \eqref{eq:GradPotencial}, $\dot p_\rho(t)=-\kappa(\mathbf q(t))\rho(t)<0,$ for $t\in I_{\mathbf z_0}^+$, because $\kappa(\mathbf q(t))>0$ and $\rho(t)>0$. Therefore, $p_\rho(t)\leq p_\rho(0)=p_{\rho_0}<0.$ Since $\dot\rho=p_\rho$, it follows that
\[
\rho(t)
=
\rho_0+\int_0^t p_\rho(s)\,\dd s
\leq
\rho_0+p_{\rho_0}t.
\]

Set $T:=-{\rho_0}/{p_{\rho_0}}>0.$ If $\omega_+>T$, then $T\in I_{\mathbf z_0}^+$ and the previous inequality gives $\rho(T)\leq \rho_0+p_{\rho_0}T=0$, which contradicts $\rho(T)>0$. Consequently, we obtain
\[
0<\omega_+\leq T
=
-\frac{\rho_0}{p_{\rho_0}}<+\infty.
\]
Since $\omega_+$ is finite, Theorem~\ref{thm:finite-time-singularities} implies that
\[
\lim_{t\rightarrow\omega_+}
\delta\bigl(\mathbf q(t)\bigr)=0.
\]
Moreover, using Theorem~\ref{theorem:approach-implies-collision} we have that there exists a point $\mathbf q_*\in\mathcal{S}$ such that
\[
\lim_{t\rightarrow\omega_+}\mathbf q(t)=\mathbf q_*.
\]
\end{proof}





\begin{remark}
The strict inequality \(p_{\rho}(0)<0\) in
Theorem~\ref{theorem:SufficientCondition} is only needed in order to
apply the theorem at the initial time. Suppose instead that $p_{\rho}(0)=0$. Since \(\rho_0>0\), the cone condition implies
\(\mathbf p_0=\mathbf 0\). Moreover, as long as \(\rho(t)>0\),
equation~\eqref{eq:GradPotencial} gives
\[
\dot p_\rho(t)=-\kappa(\mathbf q(t))\rho(t)<0.
\]
Consequently, $p_\rho(t)<0$ for every $t>0$ before the first intersection with the axis \(\rho=0\). It remains to verify that the cone condition is satisfied at positive
times. Let
\[
\chi_{\mathbf p}(t)
=
\xi(t)+\frac{p_\xi(t)}{\lambda(t)},
\qquad
\lambda(t):=-\frac{p_\rho(t)}{\rho(t)}.
\]
Using \(p_\rho(0)=p_\xi(0)=0\) and equation~\eqref{eq:GradPotencial},
we obtain
\[
p_\xi(t)
=
\kappa(\mathbf q_0)
\bigl(\chi(\mathbf q_0)-\xi_0\bigr)t+o(t),\quad 
p_\rho(t)
=
-\kappa(\mathbf q_0)\rho_0t+o(t).
\]
Therefore, $\lambda(t)=\kappa(\mathbf q_0)t+o(t)$, and hence
\[
\lim_{t\downarrow0}\chi_{\mathbf p}(t)
=
\chi(\mathbf q_0)\in I.
\]
The differential equation
\[
\dot\chi_{\mathbf p}
=
\frac{\kappa(\mathbf q)}{\lambda}
\bigl(\chi(\mathbf q)-\chi_{\mathbf p}\bigr)
\]
then shows that \(\chi_{\mathbf p}(t)\in I\) for every positive time
before the first intersection with \(\rho=0\). The variation-of-constants argument gives
\[
\mathbf q(t)\in\mathcal T(\mathbf q_0),
\qquad
\mathbf p(t)\in PWC(\mathbf q(t)).
\]
Consequently, we have $(\mathbf q(t),\mathbf p(t))\in\mathcal C(\mathbf q_0).$ Moreover, since
\(\mathbf q(t)\in\mathcal T(\mathbf q(t))\) by taking
\(\lambda_1=0\), we also have
\[
(\mathbf q(t),\mathbf p(t))\in\mathcal C(\mathbf q(t)).
\]

Thus, although Theorem~\ref{theorem:SufficientCondition} cannot be
applied at \(t=0\) when \(p_\rho(0)=0\), it can be applied to the
time-shifted solution at any time \(\tau>0\) for which the solution
remains in the upper half-plane. In particular, the remaining collision
time satisfies
\[
\omega_+-\tau
\leq
-\frac{\rho(\tau)}{p_\rho(\tau)}.
\]
\end{remark}


The following theorem states that the final arc of every trajectory ending in an interior
collision satisfies the cone condition.

\begin{theorem}[Final cone condition for interior collisions]
\label{theorem:InteriorCollisionCone}
Let $(\mathbf q(t),\mathbf p(t))=(\xi(t),\rho(t)),p_\xi(t),p_\rho(t))$ be a singular solution colliding with an interior point of the segment. That is to say
\[
\mathbf q(t)\longrightarrow
\mathbf q^*=(\xi^*,0)
\qquad\text{as }t\longrightarrow\omega_+,
\]
where $\xi^*\in(\xi_-,\xi_+).$ Then there exists \(t_*<\omega_+\) such that $\mathbf p(t)\in\mathcal C(\mathbf q(t))$ for every $t\in[t_*,\omega_+).$
\end{theorem}

\begin{proof}
Since the trajectory collides with the segment, the axial angular momentum
vanishes and the motion takes place in the invariant plane. In this plane, we restrict to the upper half-plane $\rho>0$ ($\rho<0$ is analogous) and we have
\[
\dot\rho=p_\rho,
\qquad
\dot p_\rho
=
-G\rho
\int_{\xi_-}^{\xi_+}
\frac{\sigma_{\mathrm c}(\chi)}
{\bigl((\xi-\chi)^2+\rho^2\bigr)^{3/2}}
\,d\chi<0,\qquad
\dot p_\xi
=
G
\int_{\xi_-}^{\xi_+}
\frac{\sigma_{\mathrm c}(\chi)(\chi-\xi)}
{\bigl((\xi-\chi)^2+\rho^2\bigr)^{3/2}}
\,d\chi.
\]
Thus \(p_\rho\) is strictly decreasing. Since
\(\rho(t)>0\) for \(t<\omega_+\) and \(\rho(t)\to0\), there exists
\(t_0<\omega_+\) such that $p_\rho(t)<0$ for every $t\in[t_0,\omega_+)$. For \(p_\rho<0\), define the point where the ray
\(\mathbf q+\lambda\mathbf p\), \(\lambda\geq0\), intersects the axis
\(\rho=0\) by
\[
\chi_{\mathbf p}(t)
:=
\xi(t)-\rho(t)\frac{p_\xi(t)}{p_\rho(t)}.
\]
By the definition of the attraction cone, $\mathbf p(t)\in\operatorname{PWC}(\mathbf q(t))$ is equivalent to $\xi_-\leq\chi_{\mathbf p}(t)\leq\xi_+.$ We shall prove that
\[
\chi_{\mathbf p}(t)\longrightarrow\xi^*.
\]
Hence, assuming this convergence is true, there exist $t_*<\omega_+$ such that $\chi_{\mathbf p}(t)\in(\xi_-,\xi_+)$. Finally, applying the forward-invariance result of Lemma \ref{lem:tecnico1} at such a time
\(t_*\), we conclude that
\[
\mathbf p(t)\in\operatorname{PWC}(\mathbf q(t)), \qquad t\in[t_*,\omega_+).
\]
Therefore, it only remains to show that $\chi_{\mathbf p}(t)\longrightarrow\xi^*$ to complete the proof of this theorem. This is equivalent to proving that $$\frac{p_\xi(t)}{p_\rho(t)}\rightarrow0,$$
which we establish in three steps as follows.

\noindent
\textit{Step 1:} For the estimates below, write $\xi=\xi(t),$ and $\rho=\rho(t),$ and define
\[
J(t):=
\int_{\xi_-}^{\xi_+}
\frac{\sigma_{\mathrm c}(\chi)}
{\bigl((\xi-\chi)^2+\rho^2\bigr)^{3/2}}
\,d\chi,\qquad K(t):=
\int_{\xi_-}^{\xi_+}
\frac{\sigma_{\mathrm c}(\chi)(\chi-\xi(t))}
{\bigl((\xi(t)-\chi)^2+\rho(t)^2\bigr)^{3/2}}
\,d\chi.
\]
We rescale the integration parameter in $J(t)$ by means of the following change of variable
\[
u=\frac{\chi-\xi(t)}{\rho(t)},\quad \chi=\xi(t)+\rho(t)u,\quad d\chi=\rho(t)\,du.
\]
Hence, we have $\bigl((\xi(t)-\chi)^2+\rho(t)^2\bigr)^{3/2}=\rho(t)^3(1+u^2)^{3/2}$, and considering that $\xi(t)\longrightarrow\xi^*\in(\xi_-,\xi_+)$ and $\rho(t)\longrightarrow0$, we obtain
\[
\rho(t)^2J(t)
=
\int_{\frac{\xi_- -\xi(t)}{\rho(t)}}^{\frac{\xi_+ -\xi(t)}{\rho(t)}}
\frac{\sigma_{\mathrm c}\bigl(\xi(t)+\rho(t)u\bigr)}
{(1+u^2)^{3/2}}
\,du\longrightarrow\sigma_{\mathrm c}(\xi^*)
\int_{-\infty}^{+\infty}
\frac{du}{(1+u^2)^{3/2}}=2\sigma_{\mathrm c}(\xi^*),
\]
Substituting this expression into the radial equation yields
\begin{equation}
    \label{eq:prho}
    \dot p_\rho(t)=-\frac{2G\sigma_{\mathrm c}(\xi^*)+o(1)}{\rho(t)}.
\end{equation}

We next estimate the axial force. Choose \(\delta>0\) such that $[\xi^*-2\delta,\xi^*+2\delta]\subset I$. For all sufficiently large \(t\), $[\xi(t)-\delta,\xi(t)+\delta]\subset I$. Using \(s=\chi-\xi(t)\), we can write
\[
K(t)
=
\int_{-\delta}^{\delta}
\frac{\sigma_{\mathrm c}(\xi(t)+s)s}
{(s^2+\rho(t)^2)^{3/2}}
\,ds
+
R_\delta(t),
\]
where \(R_\delta(t)\) is the contribution from
\(\lvert s\rvert\geq\delta\). Since the denominator is bounded away from
zero on this part of the integral, we have $\lvert R_\delta(t)\rvert\leq C_\delta$ for some constant \(C_\delta\) independent of \(t\). Moreover, we have
\[
\int_{-\delta}^{\delta}
\frac{\sigma_{\mathrm c}(\xi(t))s}
{(s^2+\rho(t)^2)^{3/2}}
\,ds
=0.
\]
because the integrand is odd. Therefore,
\[
K(t)
=
\int_{-\delta}^{\delta}
\frac{
\bigl[
\sigma_{\mathrm c}(\xi(t)+s)-\sigma_{\mathrm c}(\xi(t))
\bigr]s
}
{(s^2+\rho(t)^2)^{3/2}}
\,ds
+
R_\delta(t),
\]
and using the continuity of the density function, we have the following bound 
\[
\left|
\sigma_{\mathrm c}(\xi(t)+s)
-
\sigma_{\mathrm c}(\xi(t))
\right|
\leq
\omega_{\sigma_{\mathrm c}}(\delta)
\]
for \(\lvert s\rvert\leq\delta\). Hence,
\[
\begin{aligned}
\rho(t)\lvert K(t)\rvert
&\leq
\omega_{\sigma_{\mathrm c}}(\delta)
\rho(t)
\int_{-\delta}^{\delta}
\frac{\lvert s\rvert}
{(s^2+\rho(t)^2)^{3/2}}
\,ds
+
\rho(t)\lvert R_\delta(t)\rvert \leq
2\omega_{\sigma_{\mathrm c}}(\delta)+o(1),
\end{aligned}
\]
because
\[
\rho
\int_{-\delta}^{\delta}
\frac{\lvert s\rvert}{(s^2+\rho^2)^{3/2}}
\,ds
=
2\left(
1-\frac{\rho}{\sqrt{\delta^2+\rho^2}}
\right)
\leq2.
\]
Taking the limit \(t\to\omega_+\), and then letting \(\delta\to0\), we have $2\omega_{\sigma_{\mathrm c}}(\delta)\rightarrow0$, which gives
\begin{equation}
    \label{eq:r}
    r(t):=\frac{\dot p_\xi(t)}{\dot p_\rho(t)}\leq -\dfrac{2\omega_{\sigma_{\mathrm c}}(\delta)+o(1)}{2G\sigma_{\mathrm c}(\xi^*)+o(1)}
\longrightarrow0.
\end{equation}

\noindent\textit{Step 2:}
We now show that \(p_\rho(t)\longrightarrow-\infty\). Since \(p_\rho\) is
decreasing and eventually negative, it has a limit in
\([-\infty,0)\). Suppose, by contradiction, that this limit is finite $p_\rho(t)\longrightarrow p_{\rho,\infty}<0$. Since \(\rho(t)\to0\) and \(\omega_+-t\to0\), l'Hôpital's rule gives
\[
\lim_{t\to\omega_+}
\frac{\rho(t)}{\omega_+-t}
=
\lim_{t\to\omega_+}
\frac{\dot\rho(t)}{-1}
=
-p_{\rho,\infty}>0.
\]
Hence, $\rho(t)\sim (-p_{\rho,\infty})(\omega_+-t)$, and combining this relation with \eqref{eq:prho}
we obtain
\[
\dot p_\rho(t)
\sim
-\frac{C}{\omega_+-t},
\qquad
C:=
\frac{2G\sigma_{\mathrm c}(\xi^*)}
{-p_{\rho,\infty}}>0.
\]

Therefore, there exist $t_1<\omega_+$ such that
\[
\dot p_\rho(t)
\leq
-\frac{C}{2(\omega_+-t)},\qquad t\in(t_1,\omega_+).
\]
Integrating from \(t_1\) to \(t\) gives
\[
\begin{aligned}
p_\rho(t)-p_\rho(t_1)
&=
\int_{t_1}^t\dot p_\rho(s)\,ds \leq
-\frac{C}{2}
\int_{t_1}^t\frac{ds}{\omega_+-s}\rightarrow-\infty.
\end{aligned}
\]
Consequently, $p_\rho(t)\longrightarrow-\infty$, which contradicts the assumed finite limit \(p_{\rho,\infty}\).

\noindent\textit{Step 3:} For $r(t)$ defined in \eqref{eq:r} and $\epsilon>0$, we choose $t_2<\omega_+$ such that
\[
\lvert r(t)\rvert\leq\varepsilon
\qquad
\text{for all }t\in[t_2,\omega_+).
\]
Then, since \(\dot p_\rho<0\),
\[
\begin{aligned}
\left|p_\xi(t)-p_\xi(T)\right|
&=
\left|
\int_T^t r(s)\dot p_\rho(s)\,ds
\right| \leq
\varepsilon
\int_T^t-\dot p_\rho(s)\,ds =
\varepsilon\bigl(p_\rho(T)-p_\rho(t)\bigr).
\end{aligned}
\]
Dividing by \(-p_\rho(t)\) and using
\(p_\rho(t)\to-\infty\), we obtain
\[
\limsup_{t\to\omega_+}
\left|
\frac{p_\xi(t)}{p_\rho(t)}
\right|
\leq\varepsilon.
\]
Since \(\varepsilon>0\) is arbitrary,
\[
\frac{p_\xi(t)}{p_\rho(t)}
\longrightarrow0.
\]

Finally, we define $t_*:=min\{t_1,t_2\}$.
\end{proof}


\subsection{Conditional quadratic refinement of the
collision-time estimate}
\label{subsec:quadratic-bound}

The linear collision-time estimate obtained in Theorem~5.3 follows from the negativity of the transverse acceleration, $\dot p_\rho(t)=\ddot{\rho}(t)<0.$ A sharper upper bound can be obtained if the additional condition $\dot p_\rho(t)\leq \dot p_{\rho0}<0$ holds throughout the interval preceding collision, where $\dot p_{\rho0}=\dot p_\rho(0)$ denotes the initial transverse acceleration. Integrating twice under this condition yields
\[
p_\rho(t)\leq p_{\rho0}+\dot p_{\rho0}t, \qquad \rho(t)\leq
\rho_0+p_{\rho0}t+
\frac{1}{2}\dot p_{\rho0}t^2.
\]
The positive zero of the quadratic expression gives
\begin{equation}
T_{\mathrm{quad}}
=
\frac{
p_{\rho0}+
\sqrt{p_{\rho0}^{\,2}
-2\dot p_{\rho0}\rho_0}
}{
-\dot p_{\rho0}
}
\label{eq:quadratic-collision-bound}
\end{equation}

Consequently, whenever the acceleration comparison holds, the collision time satisfies
\[
\tau_*\leq T_{\mathrm{quad}}
<
-\frac{\rho_0}{p_{\rho0}}.
\]

The acceleration condition is not guaranteed for arbitrary positive density distributions. However, it can be established for a homogeneous segment when the initial conditions satisfy the sufficient angular-cone criterion, as shown below.

\begin{proposition}
\label{propo:QuadraticBound}
Consider a homogeneous segment with constant linear density $\sigma_c(\chi)\equiv\sigma_0>0$ for $\chi\in I$. Let $z_0=(q_0,p_0)$ be an initial condition in the invariant upper half-plane, where $q_0=(\xi_0,\rho_0)$, $\rho_0>0$, and $p_{\rho0}<0$.

If $(q_0,p_0)\in\mathcal{C}(q_0)$, then
\begin{equation}
\dot p_\rho(t)\leq\dot p_{\rho0}<0,
\qquad 0\leq t<\tau_*
\end{equation}
where $\tau_*$ denotes the collision time. Consequently, the quadratic estimate provides an upper bound for the collision time,
\begin{equation}
\tau_*\leq T_{\mathrm{quad}}
<
-\frac{\rho_0}{p_{\rho0}}.
\end{equation}
\end{proposition}

\begin{proof}
For constant linear density, the transverse acceleration is given by
\begin{equation}
\dot p_\rho(t)
=
-G\sigma_0\rho(t)
\int_{\xi_-}^{\xi_+}
\frac{d\chi}
{\left[(\xi(t)-\chi)^2+\rho(t)^2\right]^{3/2}}\nonumber
\end{equation}
Following Theorem~\ref{theorem:InteriorCollisionCone}, we introduce the change of variables $u=(\chi-\xi(t))/\rho(t)$ giving
\begin{equation}
-\dot p_\rho(t)=\frac{G\sigma_0}{\rho(t)}
\int_{u_-(t)}^{u_+(t)}
\frac{du}{(1+u^2)^{3/2}}\nonumber
\end{equation}
with $u_-(t)=({\xi_- -\xi(t)})/{\rho(t)},$ and $u_+(t)=({\xi_+ -\xi(t)})/{\rho(t)}.$ By Lemma \ref{lem:tecnico1}, the angular cone is forward invariant, and the trajectory remains within the triangle determined by $q_0$ and the endpoints of the segment. Defining
\[
\alpha(t)=\frac{\rho(t)}{\rho_0},
\qquad 0<\alpha(t)\leq1,
\]
the triangle invariance property implies
\begin{equation}
\alpha(t)\xi_0+(1-\alpha(t))\xi_-
\leq \xi(t)
\leq
\alpha(t)\xi_0+(1-\alpha(t))\xi_+.
\end{equation}
Using $\rho(t)=\alpha(t)\rho_0$, these inequalities give
\begin{equation}
u_-(t)\leq u_-(0),
\qquad
u_+(t)\geq u_+(0).
\end{equation}
Therefore, $[u_-(0),u_+(0)] \subseteq [u_-(t),u_+(t)].$ Since the integrand $(1+u^2)^{-3/2}$ is strictly positive, it follows that
\begin{equation}
\int_{u_-(t)}^{u_+(t)}
\frac{du}{(1+u^2)^{3/2}}
\geq
\int_{u_-(0)}^{u_+(0)}
\frac{du}{(1+u^2)^{3/2}}.
\end{equation}
Furthermore, the inward motion guarantees $0<\rho(t)\leq\rho_0$, and hence $1/\rho(t)\geq1/\rho_0$. Combining these inequalities yields
\begin{equation}
\begin{aligned}
-\dot p_\rho(t)
&=
\frac{G\sigma_0}{\rho(t)}
\int_{u_-(t)}^{u_+(t)}
\frac{du}{(1+u^2)^{3/2}}\geq
\frac{G\sigma_0}{\rho_0}
\int_{u_-(0)}^{u_+(0)}
\frac{du}{(1+u^2)^{3/2}}=-\dot p_{\rho0}.
\end{aligned}
\end{equation}
Thus, $\dot p_\rho(t)\leq\dot p_{\rho0}<0$ throughout the pre-collision interval. The required acceleration condition is therefore satisfied, and the quadratic collision-time estimate follows by integration, giving
\[
\tau_*\leq T_{\mathrm{quad}}
<
-\frac{\rho_0}{p_{\rho0}}.
\]
\end{proof}

\section{Numerical experiments}
\label{sec:numerical-experiments}

We now present numerical simulations illustrating the theoretical results.
The integrations use four density profiles: uniform, linear, sinusoidal, and
center-heavy density. The coordinates are expressed in the same
nondimensional units used throughout the manuscript, with lengths normalized
by the segment half-length \(L\).

Since the potential is singular on the segment, the integrations associated
with near-collision trajectories are stopped when the distance to the segment
reaches a prescribed tolerance \(\delta_{\mathrm{tol}}\). Accordingly, the
reported collision times should be interpreted as numerical near-collision
times rather than exact collision times.

The initial conditions used in the simulations are summarized in Table~\ref{tab:numerical-initial-conditions}. Cases (a)--(c) correspond to initial conditions inside the angular cone. Case (a) starts in the interior of the cone, case (b) starts on its boundary, and case (c) starts outside the orthogonal projection of the segment while remaining inside the cone. Case (d) is a planar trajectory outside the angular cone, whereas case (e) has
nonzero axial angular momentum.

\begin{figure}[htbp]
    \centering
    \includegraphics[width=17.0cm]{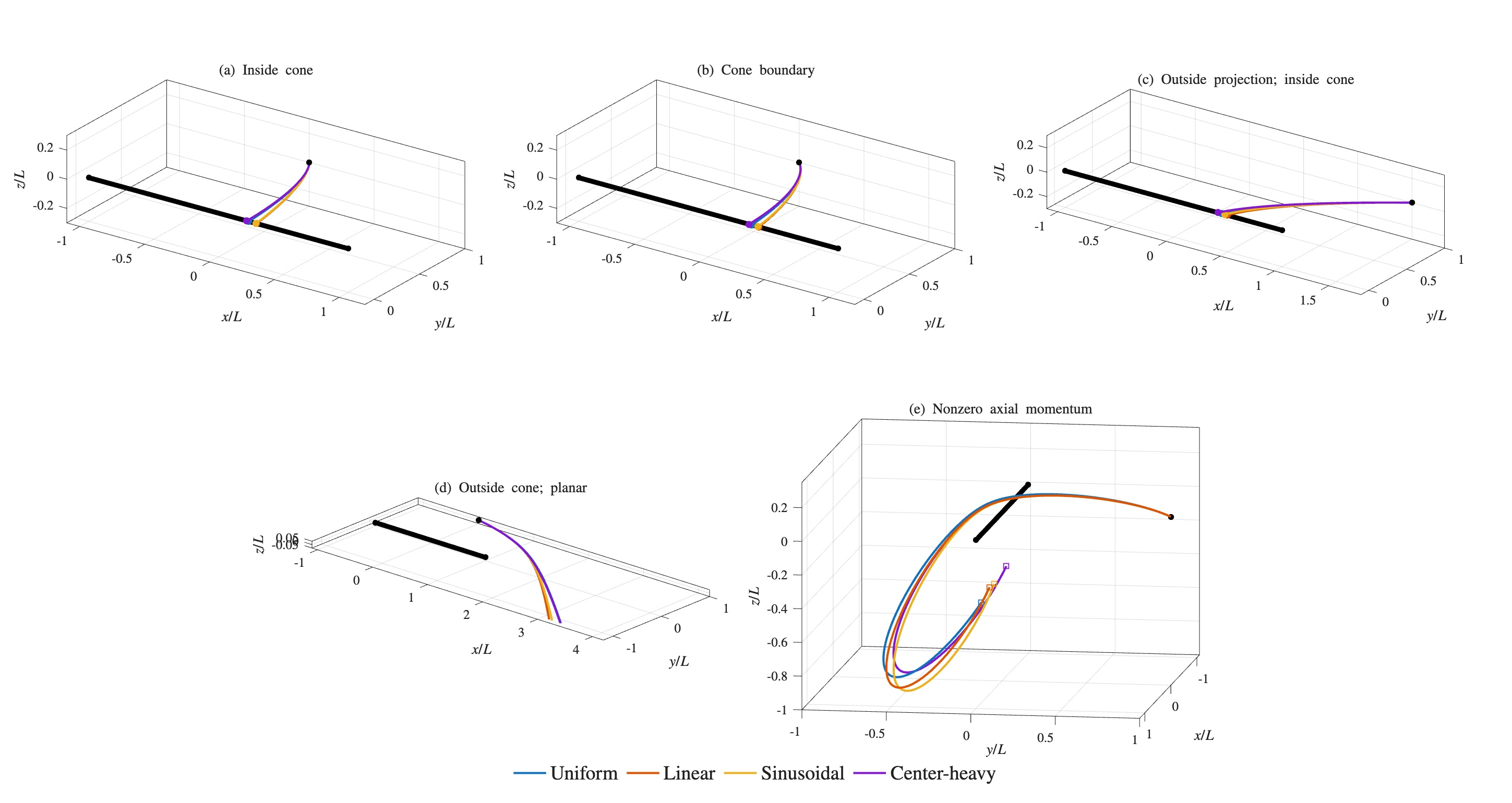}
    \caption{Numerical trajectories for the four density profiles. Panels
    (a)--(c) correspond to initial conditions with zero axial angular
    momentum and satisfying the angular-cone condition. Panel (d) shows a
    planar trajectory outside the angular cone, and panel (e) shows a
    trajectory with nonzero axial angular momentum. The blue, orange, yellow,
    and purple curves correspond to the uniform, linear, sinusoidal, and
    center-heavy densities, respectively. The integrations for panels (d)
    and (e) were continued until \(t=12\); shorter time intervals are shown
    for visualization purposes.}
    \label{fig:numerical-collisions}
\end{figure}

\begin{table}[htbp]
\centering
\caption{Initial conditions used in the numerical experiments displayed in
Figure~\ref{fig:numerical-collisions}.}
\label{tab:numerical-initial-conditions}
\begin{tabular}{c c c c p{0.34\textwidth}}
\hline
Case & \(\mathbf r_0\) & \(\mathbf v_0\) & \(\ell\)
& Theoretical prediction \\
\hline
(a)
& \((0,1,0)\)
& \((0.35,-0.50,0)\)
& \(0\)
& Collision, \(t^*\leq 2\) \\

(b)
& \((0,1,0)\)
& \((0.50,-0.50,0)\)
& \(0\)
& Collision, \(t^*\leq 2\) \\

(c)
& \((1.5,1,0)\)
& \((-0.75,-0.50,0)\)
& \(0\)
& Collision, \(t^*\leq 2\) \\

(d)
& \((0,1,0)\)
& \((2.50,-0.50,0)\)
& \(0\)
& The cone condition is not satisfied; no collision is detected for
\(0\leq t\leq 12\) \\

(e)
& \((0,1,0)\)
& \((0.35,-0.50,0.30)\)
& \(0.30\)
& Collision excluded by Theorem~\ref{thm:necessary-angular-momentum} \\
\hline
\end{tabular}
\end{table}

For the numerical simulations in cases (a)--(c), the initial radial coordinate and radial momentum are $\rho_0=1$ and $p_{\rho0}=-0.50.$ Therefore, the collision-time estimate from
Theorem~\ref{theorem:SufficientCondition} is
\[
-\frac{\rho_0}{p_{\rho0}}=2.
\]
The corresponding numerical near-collision times are listed in
Table~\ref{tab:near-collision-times}.

\begin{table}[htbp]
\centering
\caption{Numerical near-collision times for the four density profiles.}
\label{tab:near-collision-times}
\begin{tabular}{c c c c c}
\hline
Case & Uniform & Linear & Sinusoidal & Center-heavy \\
\hline
(a) & \(0.74347\) & \(0.73500\) & \(0.73226\) & \(0.72705\) \\
(b) & \(0.74442\) & \(0.73294\) & \(0.72993\) & \(0.72933\) \\
(c) & \(0.93309\) & \(0.88581\) & \(0.89755\) & \(0.94473\) \\
\hline
\end{tabular}
\end{table}

In all three cases, the numerical near-collision times are below the linear bound established in Theorem \ref{theorem:SufficientCondition}, in agreement with the sufficient collision condition. For uniform density, the quadratic refinement derived in Section \ref{subsec:quadratic-bound} provides a sharper upper bound, which is further examined through numerical comparisons of collision times and transverse trajectories.

Case (d) illustrates that the angular momentum condition is necessary but not sufficient. The initial condition is planar having zero axial angular momentum, but it does not belong to the angular cone. No collision is detected for any of the four density profiles during the integration interval \(0\leq t\leq 12\), and the trajectories move away from the segment. 

Case (e) illustrates the role of axial angular momentum. Since
\(\ell=0.30\), Theorem~\ref{thm:necessary-angular-momentum} excludes
collision with the segment for every density profile. The corresponding
trajectories remain non-collisional throughout the numerical integration.


To further examine the quadratic refinement established in Section \ref{subsec:quadratic-bound}, we consider a uniform-density segment with $L=1$ and $G\sigma_0=1$.
A set of 1000 initial conditions was randomly sampled, together with three reference cases, all satisfying the sufficient angular-cone criterion. For each trajectory, the numerical near-collision time $T_{\mathrm{num}}$ is recorded when the transverse distance reaches $\varepsilon=10^{-6}$. The resulting times are compared with the quadratic estimates computed from the corresponding initial conditions.

\begin{figure}[hbt!]
    \centering
    \includegraphics[width=0.5\linewidth]
    {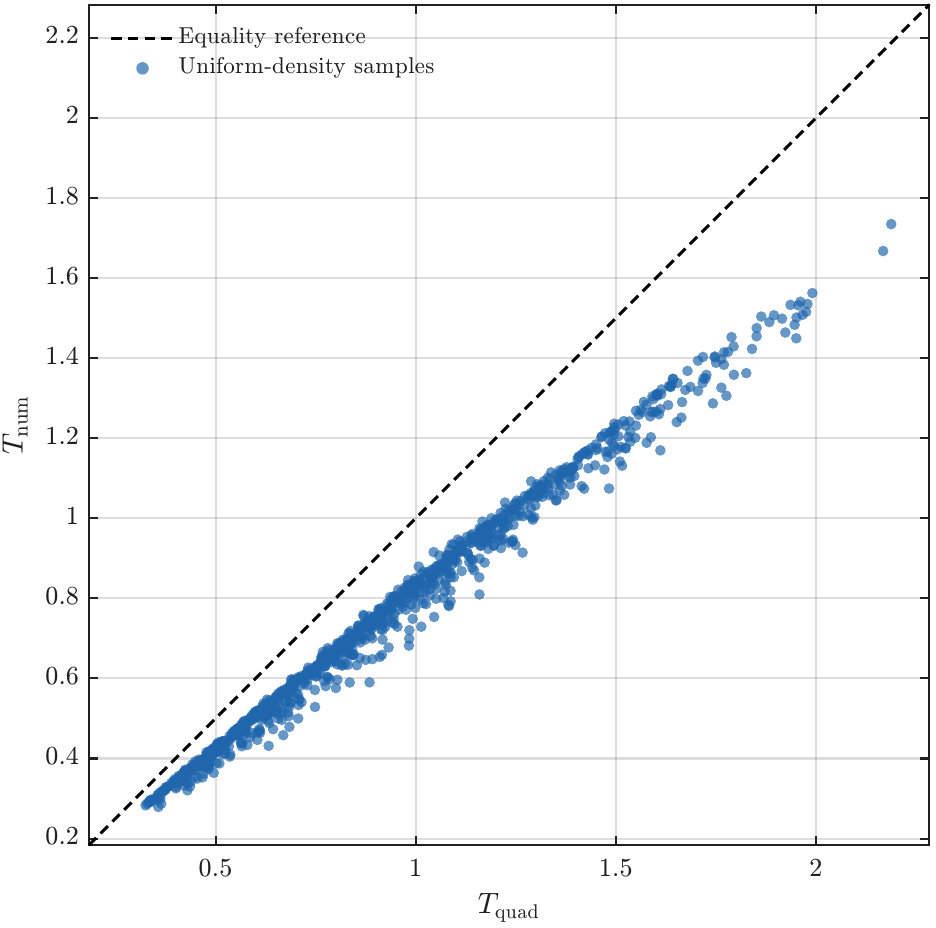}
    \caption{
    Comparison of numerical collision times
    $T_{\mathrm{num}}$ and quadratic estimates
    $T_{\mathrm{quad}}$ for initial conditions satisfying
    the sufficient angular-cone criterion under uniform density.
    The dashed diagonal represents equality.
    }
    \label{fig:quadratic_sweep}
\end{figure}

Figure~\ref{fig:quadratic_sweep} compares $T_{\mathrm{num}}$ and $T_{\mathrm{quad}}$ across the sampled trajectories. All successfully integrated samples lie below the equality line, consistent with the quadratic upper bound established analytically for uniform density.

\begin{figure}[hbt!]
    \centering
    \includegraphics[width=\linewidth]
    {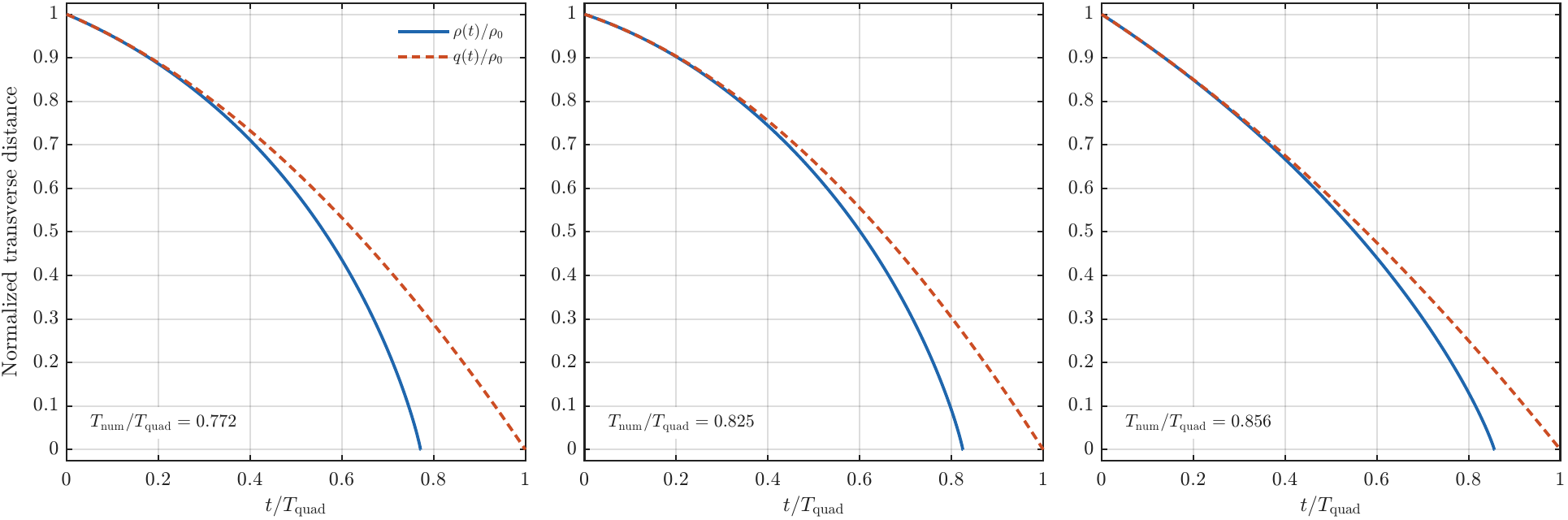}
    \caption{
    Evolution of the numerical transverse distance
    $\rho(t)$ (solid curves) and the quadratic comparison
    function $\rho_{\mathrm{quad}}(t)$ (dashed curves)
    for three representative trajectories under uniform
    density. The initial conditions satisfy the sufficient
    angular-cone criterion and correspond approximately
    to the 10th, 50th, and 90th percentiles of
    $T_{\mathrm{num}}/T_{\mathrm{quad}}$.
    Time and transverse distance are normalized by
    $T_{\mathrm{quad}}$ and $\rho_0$, respectively.
    The numerical trajectories terminate at $\rho=10^{-6}$, while the quadratic comparison functions reach zero at $t/T_{\mathrm{quad}}=1$.
    }
    \label{fig:uniform-quadratic-trajectories}
\end{figure}

To examine the comparison at the trajectory level, Figure~\ref{fig:uniform-quadratic-trajectories} shows the evolution of the transverse distance for three representative initial conditions. These are selected from the successfully integrated samples according to the 10th, 50th, and 90th percentiles of the ratio $T_{\mathrm{num}}/T_{\mathrm{quad}}$.

For each trajectory, the numerical solution $\rho(t)$ is compared with the quadratic function
\[
\rho_{\mathrm{quad}}(t)
=
\rho_0+p_{\rho0}t
+\frac{1}{2}\dot p_{\rho0}t^2.
\]
The numerical transverse distance remains below the quadratic comparison function throughout the computed pre-collision interval, illustrating the inequality underlying the collision-time estimate. The separation between the curves varies among the three cases, reflecting different degrees of conservativeness of the quadratic bound. In each case, the numerical collision occurs before the quadratic comparison function reaches zero.

\FloatBarrier



\section*{Acknowledgment}
J.~V. was partially supported by ANID--FONDECYT Iniciación Grant No.~11240582. J.~L.~Z. was supported by ANID--FONDECYT Iniciación Grant No.~11250703.

\bibliographystyle{plain}
\bibliography{Bibliography}

\end{document}